\newtheorem{theorem}{\sc Theorem}
 
\newtheorem{proposition}[theorem]{\sc Proposition}
\newtheorem{lemma}[theorem]{\sc Lemma}
\newtheorem{corollary}[theorem]{\sc Corollary}

\newcommand{\proof}{\noindent {\sc Proof. }} 
 
\def\rest{\mathord{\restriction}}
\newcommand{\force}{\Vdash}
\newcommand{\open}{\Bbb}
 
\newlength{\labparwidth}
\newcommand{\cf}{\mbox{\rm cf}}

\renewcommand{\hom}{\mbox{\rm Hom}}
\newcommand{\ext}{\mbox{\rm Ext}}
\newcommand{\se}{\subseteq}
\newcommand{\set}[2]{\{#1 \colon #2\}} 
\newcommand{\fin}{$\Box$\par\medskip} 
\newcommand{\dmd}{\diamondsuit}
\renewcommand{\to}{\rightarrow}

\newcommand{\ga}{\alpha}
\newcommand{\gb}{\beta}
\newcommand{\grg}{\gamma}
\newcommand{\gd}{\delta}

\newcommand{\gk}{\kappa}
\newcommand{\gl}{\lambda}

\newcommand{\gt}{\tau}

\newcommand{\gf}{\varphi}

\newcommand{\go}{\omega}
 
\newcommand{\ha}{\aleph}
 
\newcommand{\oP}{{\open P}}
\newcommand{\oQ}{{\open Q}}
\newcommand{\oR}{{\open R}}
\newcommand{\oZ}{{\open Z}}
 
\documentstyle[12pt,amsfonts,amssymb,titlepage]{article}

\title{Every Coseparable Group May Be Free}
\author{Alan H. Mekler\thanks{Research supported by NSERC. Research on
this paper was begun while the first author was visiting the Hebrew
University.}\\Department of Mathematics and Statistics\\ Simon Fraser
University\\Burnaby, B.C.,V5A 1S6\\ Canada \and Saharon
Shelah\thanks{Publication
\#418. Research supported by the BSF.}\\Institute of Mathematics\\The
Hebrew University\\Jerusalem, Israel\\and\\Department of Mathematics\\ Rutgers
University\\New Brunswick, New Jersey 08903\\USA} 
\date{}

\newcommand{\Fn}{{\rm Fn}}
\newcommand{\tA}{\tilde{A}}
\newcommand{\tC}{\tilde{C}}
\newcommand{\tQ}{\tilde{Q}}
\newcommand{\tS}{\tilde{S}}
\newcommand{\tf}{\tilde{f}}
\newcommand{\tT}{\tilde{T}}
\newcommand{\ax}{${\rm Ax}^+(\ha_1\mbox{-complete forcing})$}
\begin{document}
\maketitle

\begin{abstract}
We show that if $2^{\ha_0}$ Cohen reals are added to the universe,
then for every reduced non-free torsion-free abelian group $A$ of
cardinality less than the continuum, there is a prime $p$ so that
$\ext_p(A, {\open Z}) \neq 0$. In particular if it is consistent that
there is a supercompact cardinal, then it is consistent (even with
weak CH) that every coseparable group is free. The use of some large
cardinal hypothesis is needed.
\end{abstract}

\section{Introduction}

	There are several approximations to being free for abelian
groups.  (From here on by {\em group} we will mean {\em abelian
group}.)  For some of these notions such as being hereditarily
separable or being a Whitehead group it is well known that it is
independent from the usual axioms of set theory whether or not every
such group is free.  (A group is {\em separable} if every finite set
is contained in a free subgroup which is a direct summand and {\em
hereditarily separable} if every subgroup is separable. A group $A$ is
a {\em Whitehead group} if $\ext(A, \oZ) = 0 $. Any Whitehead group is
hereditarily separable.)  On the other hand it can be shown in ZFC
that there are non-free groups of cardinality $\ha_1$ which are
$\ha_1$-separable (i.e., any countable subset is contained in free
countable direct summand). A somewhat mysterious class is the class
of coseparable groups (to be defined later). Any hereditarily separable
group is coseparable while every coseparable group is separable. So
this class lies between those which exhibit independence phenomena and
those for which we have absolute existence results.

	Early on Chase \cite{Ch} showed that CH implied the existence
of a non-free coseparable group of cardinality $\ha_1$.  Later
Sageev-Shelah \cite{SS} and Eklof-Huber \cite{EH} elaborated this
result to construct non-free groups of cardinality $\ha_1$ with
specified $p$-rank of Ext (as $p$ ranges over the primes).  Their
proofs continued to use CH\@.  This reliance on CH was in itself
unusual.  In the cases where the  existence of a group was known to be
independent, the existence was also independent of CH\@.  As well many
results in abelian group theory that follow from CH can also be shown
to follow from weak CH, i.e., $2^{\ha_0} < 2^{\ha_1}$.  In this paper
we will show that the use of CH is necessary.  Namely we will show
that is consistent with weak CH (assuming as always the consistency of
ZFC) that every coseparable group of cardinality $\ha_1$ is free.
Furthermore if it is consistent that a supercompact cardinal exists
then it is consistent (with weak CH if desired) that every coseparable
group is free.  A reference for the facts mentioned above and the ones
we will quote below is \cite{EM}.

	A group is {\em $\ha_1$-free} if every countable subgroup is
free and, more generally, for an uncountable cardinal $\gk$ a group is
{\em $\gk$-free} if every subgroup of cardinality less than $\gk$ is
free. A group $A$ is {\em coseparable} if it is $\ha_1$-free and every
subgroup $B$ such that 
$A/B$ is finitely generated contains a direct summand $C$ of $A$ so
that $A/C$ is finitely generated. There are several equivalents to
being coseparable. The one we will use throughout is that a group $A$
is coseparable if and only if $\ext(A, \oZ)$ is 
torsion-free or equivalently $\ext_p(A, \oZ) = 0$ for all primes $p$.
Here $\ext_p(A, \oZ)$ denotes the subgroup of $\ext(A, \oZ)$
consisting of those elements whose order is a power of $p$. There is a
useful characterization of when $\ext_p(A, \oZ) = 0$. For a group $A$
and a prime $p$, $\ext_p(A, \oZ) = 0$ if and only if for every
homomorphism $h\colon A \to \oZ/p\oZ$ there is a homomorphism
$\hat{h}\colon A \to \oZ$ (called a {\em lifting} of $h$) so that
$\hat{h}/p = h$. (If $f$ is a homomorphism to $\oZ$, then $f/p$
denotes the composition of $f$ with the natural map to $\oZ/p\oZ$.)

	Although (at least in some models of set theory) a subgroup of
a coseparable group is not necessarily coseparable, every pure subgroup
of a coseparable group is coseparable. The following proposition sums
up the discussion.

\begin{proposition}
\label{observation}
\begin{enumerate}
\item If $\gl$ is the minimal cardinality of a non-free coseparable group
then every coseparable group is $\gl$-free.
\item Any homomorphism from a subgroup $H$ of $K$ to $\oZ/p\oZ$ can be
extended to a homomorphism from $K$ to $\oZ/p\oZ$ provided that $H$ is
$p$-pure in $K$, i.e., for any $a\in H$, $a \in pH$ if and only if $a
\in pK$.
\end{enumerate}
\end{proposition}

\section{Proofs}
Most of our set theoretic notation will be standard. If $\gk$ is an
infinite cardinal  we will use $\Fn(I, J, \gk)$ to denote the set of partial
functions from $I$ to $J$ whose domain have cardinality less than
$\gk$. The poset $\Fn(\gk, 2, \go)$ is usually called the poset (or
forcing) for adding $\gk$ Cohen reals. As is customary we will say
that (at least)  $\gk$ Cohen reals are added to a set theoretic
universe if we force with $\Fn(\gk, 2, \go)$ (or $\Fn(\mu, 2, \go)$
for some $\mu \geq \gk$)
The principal lemma we will use is the following.

\begin{lemma}
\label{main}
Suppose $\gk$ is a cardinal and $\oP = \Fn(\gk, 2, \go)$. Then $\oP$
forces that every coseparable group of cardinality less than $\gk$ is
free.
\end{lemma}

	This lemma has as immediate consequences our main theorems.

\begin{theorem}
It is consistent with both $2^{\ha_0} < 2^{\ha_1}$ and $2^{\ha_0} =
2^{\ha_1}$, that every coseparable group of cardinality $\ha_1$ is free.
\end{theorem}

\begin{theorem}
\label{weak}
If it is consistent that a supercompact cardinal exists then it is consistent
with both $2^{\ha_0} < 2^{\ha_1}$ and $2^{\ha_0} =
2^{\ha_1}$, that every coseparable group  is free.
\end{theorem}

\proof Shelah (see~\cite{BD}) has shown that if $\gk$ is a supercompact
cardinal and $\gk$ Cohen reals are added to the universe then every
$\gk$-free group of cardinality $\gk$ is free. So after adding $\gk$
Cohen reals we have, by Lemma~\ref{main} and
Proposition~\ref{observation}, that every coseparable group is
$\gk$-free and hence free. As well we have that $2^{\ha_0} = \gk =
2^{\ha_1}$.

	For the other part of the theorem, it suffices to know that the
consistency of a supercompact cardinal implies the consistency of a
cardinal $\gk$ such that $\gk < 2^{\ha_1}$ and if $\gk$ Cohen reals are
added to the universe then every $\gk$-free group is free.  This result
is proved  in Corollary~\ref{free}. \fin

	We will also show that the use of some large cardinal is
necessary. (The consistency of the statement ``every coseparable group
is free'' implies the consistency of the existence of measurable
cardinals.) Also we will show by a more complicated forcing that
assuming the consistency of a supercompact cardinal it is consistent
that every coseparable group is free and the continuum is as small as
possible, i.e., $\ha_2$.

	There are various cases to consider in the proof of the main
lemma.  At various times we will have to consider the case of a free
group $F$ and its $p$-adic completion.  Fix a set of free generators $X$
for $F$.  Then any element $y$ of the $p$-adic completion can be written
uniquely as $\sum p^n y_n$ where each $y_n$ is a linear combination of
elements of $X$ such that all the coefficients are between $0$ and $p -
1$ (inclusive).  If we have fixed a set of free generators then given
$y$, $y_n$ will always denote the element above.  First we prove two
easy and useful propositions.

\begin{proposition}
\label{1}
Suppose $A$ is contained in the $p$-adic completion of a free group $F$
and $X$ is a set of free generators for $F$. Assume $h \in \hom(A,
{\open Z}/p{\open Z})$ and $\hat{h}$ is a lifting to $\open Z$. If $x,
y \in A$,  $\hat{h}(x) = \hat{h}(y)$ and $n$ is the least element so
that $x_n \neq y_n$ then $h(x_n) = h(y_n)$.
\end{proposition}

\proof By subtracting $\sum_{i<n} p^ix_i$ from each side and then
dividing by $p^n$, we have $\hat{h}(\sum_{n\leq m < \infty} p^{m - n}x_m) =
\hat{h}(\sum_{n\leq m<\infty} p^{m - n}y_m)$. So $\hat{h}(x_n) \equiv
\hat{h}(y_n) \bmod p$ \fin

\begin{proposition}
\label{abs-prop}
Let $A$, $F$ and $X$ be as above. Assume that there is a finite set $Y$
of elements which are not divisible by $p$ and a set $B \se A$ so that
$|B| > |F|$ and for all $x \neq y \in B$, $x_n - y_n \in Y$, where $n$
is the least natural number so that $x_n
\neq y_n$. Then $\ext_p(A, {\open Z}) \neq 0$.
\end{proposition}

\proof Suppose not. For each  $a \in Y$, choose $h_a$ so
that $h_a(a) \neq 0$. Let $\hat{h}_a$ be a lifting of $h_a$. Choose
$B' \se B$ so that $|B'| > |F|$ and for all $x, y \in B'$ and $a \in
Y$, $\hat{h}_a(x) = \hat{h}_a(y)$. By the last proposition for $x \neq
y \in B'$ and $a \in Y$, $x_n - y_n \neq a$, where $n$ is the least
natural number so that $x_n \neq y_n$. This contradicts the choice of
$Y$.
\fin

	 We will have occasion to use the following  {\em ad
hoc} definition. If $A$ is an abelian group a function $h\colon
A \to \oZ/p\oZ$ is an {\em counterexample to the vanishing of\/
$\ext_p(A, \oZ)$} if $h$ has no lifting to a map to $\oZ$.

\begin{theorem}
\label{dens-thm}
Suppose $A$ is contained in the $p$-adic completion of a free group $F$
and $|A| > |F|$. Then after adding at least $|F|$ Cohen reals to the universe
$\ext_p(A, {\open Z}) \neq 0$. 
\end{theorem}

\proof By splitting the forcing into two parts (if necessary) we can
assume that exactly $|F|$ Cohen reals are added. Let $X$ be a set of
free generators of $F$ and let $h$ be a Cohen generic function from
$X$ to $p$.  More exactly, we can assume the Cohen reals are added by
forcing with $\Fn(X, p, \go)$. Let $h$ be a generic set for this
forcing and we will also let $h$ denote the unique extension to a
homomorphism from $A$ to $\oZ/p\oZ$.  (There are many other schemes to
produce a generic function.) Suppose $\hat{h}$ is forced (by some
condition) to be a lifting of $h$. For each $a \in A$, choose a
condition $q_a$ and an integer $m_a$ so that $q_a
\force \hat{h}(a) = m_a$.  Since $|A| > \max\{|X|, |\oZ|\}$, there is
a condition $q$ so
that for some $m$ and some $B$ of cardinality greater than $|F|$, $q \force
\hat{h}(x) = m$, for all $x \in B$.  Let $X_0$ be the domain of $q$ and let
$Y$ be the set of linear combinations of elements of $X_0$ whose
coefficients are between $-p + 1$ and $p-1$. 

	Suppose first that for all $x, y \in B$, $x_n - y_n \in Y$
where $n$ is the least natural number such that $x_n \neq y_n$. Then
we are done by Proposition~\ref{abs-prop}. So we can choose $x, y \in B$ so
that this is not the case. Choose now $q_1$ extending $q$ so that
$q_1 \force h(x_n - y_n) \neq 0$. This contradicts Proposition~\ref{1}. So no
such lifting exists. \fin 

	The non-existence of a group $A$ which is coseparable and
satisfies the hypotheses of the theorem is not a consequence of
ZFC. It is consistent that there is a Whitehead group of cardinality
$\ha_1$ which is contained in the 2-adic completion of a countable set
\cite{Sh}. In fact the non-free coseparable group constructed by Chase
\cite{Ch}, assuming CH, is contained in the $\oZ$-adic completion of a
countable set.

	We can strengthen the the previous theorem to allow us to
calculate the $p$-rank.

\begin{theorem}
\label{stdens-thm}
Suppose $A$ is contained in the $p$-adic completion of a free group $F$ and
$|A| > |F|$. Then if $\gl \geq |F|$ and $\gl$ Cohen reals are added 
to the universe, $|\ext_p(A, {\open Z})| \geq \gl$.
\end{theorem}

\proof There are two cases to consider. First assume that we have a
collection $\set{Y_n}{n< \go}$ of disjoint finite subsets of $F$ so that
for all $n$, no element of $Y_n$ is divisible by $p$ modulo the subgroup
generated by $\bigcup_{m\neq n} Y_m$, the subgroup generated by $\bigcup_{n
< \go} Y_n$ is $p$-pure in $F$ and for all $n$ there is a set $B$ of
cardinality greater than the cardinality of $F$ so that if $x, y \in B$
then $x_i - y_i \in Y_n$ where $i$ is the least natural number so that $x_i
\neq y_i$. In this case by (the proof of) Proposition~\ref{abs-prop} we can
find elements $a_n \in Y_n$ such that if $h\colon A \to \oZ/p\oZ$ is a
homomorphism which is non-zero on $a_n$ then $h$ does not lift.  There
are $2^{\ha_0}$ homomorphisms $\set{h_i\colon A \to \oZ/p\oZ}{i <
2^{\ha_0}}$ such that for all $i \neq j$ there is $n$ such that $h_i
- h_j(a_n) \neq 0$. So for any $i \neq j$, $h_i - h_j$ does not lift
to a map to $\oZ$ and hence does not represent 0 in $\ext_p(A, \oZ)$. As $h_i$
and $h_j$ represent different elements of $\ext_p(A, \oZ)$,
$|\ext_p(A, \oZ)| = 2^{\ha_0}$. 

	For the second case we assume that the first does not hold. In this
case we can find a pure finitely generated subgroup $Y$ of $F$ so that for
all finite subsets $Z$ of $F$ if every element of $Z$ is not divisible by $p$
modulo $Y$ then there is no set $B$ of cardinality greater than $|F|$ so
that for all $x, y \in B$, $x_n - y_n \in Z$ where $n$ is the least natural
number so that $x_n \neq y_n$. We consider the Cohen reals as giving $\gl$
generic functions $\set{h_\ga}{\ga < \gl}$ from $F$ to $\oZ/p\oZ$ which are
$0$ on $Y$. In this case 
we can repeat the second part of the argument in Theorem~\ref{dens-thm},
since for all $\ga \neq \gb$, $h_\ga - h_\gb$ is generic.  \fin

	We now turn to the proof of the main Lemma.  Since if $A$ is a
non-free coseparable group of minimal cardinality, $A$ is $|A|$-free,
it is enough to show the result for all $\gl$-free groups $A$ of
cardinality $\gl < 2^{\ha_0}$. By Shelah's singular compactness
theorem \cite{sing}, $\gl$ is regular.  We will
consider a $\gl$-filtration of $A$ (i.e., write $A$ as a union of a
continuous chain $(A_\ga\colon \ga < \gl)$ of pure subgroups of cardinality
less than $\gl$ such that for all $\ga$, if $A/A_\ga$ is not $\gl$-free
then $A_{\ga+1}/A_\ga$ is not free). There will be two cases to 
consider.  The first and easier case is dealt with by the following
Theorem. 

\begin{theorem} 
\label{nondivis}
Suppose that $|A| = \gl$ and $A$ is $\gl$-free. Also
suppose that $(A_\ga \colon \ga < \gl)$ is a $\gl$-filtration
of $A$ and for a stationary set $E$, after adding at least $\gl$ Cohen
reals, $\ext_p(A_{\ga+1}/A_\ga, {\open Z}) 
\neq 0$ for all $\ga \in E$. Furthermore assume that for all $\ga \in E$ the
non-vanishing of $\ext_p(A_{\ga+1}/A_\ga, \oZ)$ is witnessed by a
function in the ground model.  Then if we add at least $\gl$ Cohen
reals, $\ext_p(A,
\oZ)$ has rank $2^\gl$. 

	Also if weak diamond of $E$ holds then $\ext_p(A, \oZ)$ has
rank $2^\gl$, provided that for all $\ga \in E$,
$\ext_p(A_{\ga+1}/A_\ga, {\open Z})  \neq 0$.
\end{theorem}

\proof This is a fairly routine weak diamond proof.  In
the case where we add Cohen reals the proof can be done using the
definable weak diamond (see \cite{MS}).  In order to pursue the main
case we will give the proof using the Cohen reals directly.  As well
we will only show that $\ext_p(A, {\open Z}) \neq 0$. (The proof that
$\ext_p(A, \oZ)$ has rank $2^\gl$ is  similar.)  Choose a basis $X =
\cup_{\ga<\gl} X_\ga$ of $A/pA$ so that (under the obvious
identifications) $X_\ga$ is a basis of
$(A_{\ga+1}/A_\ga)/p(A_{\ga+1}/A_\ga)$. For $\ga \in E$, let
$h_{1\ga}$ be a map, in the ground model, from $X_\ga$ to ${\open
Z}/p{\open Z}$ which does not lift to a homomorphism from
$A_{\ga+1}/A_\ga$ to $\open Z$ (even in the forcing extension).  Let
$h_{0\ga}$ denote the homomorphism which is constantly $0$ on $X_\ga$.
The important point to notice here is that if $f$ is any homomorphism
from $A_\ga$ to $\open Z$ then at most one of $f/p + h_{0\ga}$ and
$f/p + h_{1\ga}$ lifts to a homomorphism of $A_{\ga+1}$ which agrees
with $f$ on $A_\ga$.  Otherwise if $f_0, f_1$ were two such liftings
then $f_1 - f_0$ would contradict the choice of $h_{1\ga}$. For $\ga
\notin E$, let $h_{0\ga}$ and $h_{1\ga}$ be any functions from $X_\ga$
to $\oZ/p\oZ$.

	Now we can use $\gl$ Cohen reals  to define a map which cannot
be lifted. We can assume that the forcing is $\Fn(\gl, 2, \go)$ and
$G$ is a generic function from $\gl$ to 2.  Let $h$ from $X$ to
${\open Z}/p{\open Z}$ be $\sum_{\ga < \gl} h_{G(\ga)\ga}$.  Suppose
$\hat{h}$ is forced to be a lifting of $h$ (without loss of generality
we can assume that it is forced to be a lifting by the empty
condition).  Since Cohen forcing is c.c.c.\  we can find $\ga \in E$
so that $\hat h \rest A_\ga$ is decided by $G\rest \ga$.  So there is
$l = 0, 1$ and a condition, $q$, whose domain is contained in $\ga$
such that $q$ forces that $\hat h \rest A_\ga$ does not extend to a
function which lifts $h \rest A_\ga + h_{l\ga}$. If we choose $r$ so
that $r$ extends $q$ and $r(\ga) = l$ then $r$ forces that $\hat h$ is
not a lifting of $h$. This is a contradiction.
 \fin


	The last case of the main lemma is the following.

\begin{theorem}
\label{divis}
Suppose that $\gl$ is a regular cardinal and $A$ is a group of
cardinality $\gl$ so that $A$ is $\gl$-free and there is a
$\gl$-filtration of $A$ so that for a stationary set $E$ and all $\ga
\in E$ there is $a \in A_{\ga+1} \setminus A_\ga$ so that $a$ is in the
$p$-adic completion of $A_\ga$. Then if we add at least $\gl$ Cohen reals
to the universe, $\ext_p(A, \oZ) \neq 0$. In fact, $|\ext_p(A, \oZ)| \geq
\gl$. 
\end{theorem} 

\proof  In view of the Theorem~\ref{dens-thm} (or Theorem~\ref{stdens-thm}
for the stronger result) we can assume that the $p$-adic closure in
$A$ of any set of cardinality $< \gl$ has size $<\gl$.  Applying
Fodor's lemma we get that all but a non-stationary subset of $E$
consists of ordinals of cofinality $\go$. (Otherwise there would be
some $\ga$ so that the $p$-adic closure in $A$ of $A_\ga$ has
cardinality $\gl$.)  Also by Theorem~\ref{dens-thm}, we can choose the
filtration so that for all $\ga$ the $p$-adic closure in $A$ of
$A_\ga$ is contained in $A_{\ga+1}$.

	Fix a basis $X$ for $A/pA$. We now redefine the filtration
of $A$.  We can assume $A$ is the union of a $\gl$-filtration where
$A_\ga = N_\ga
\cap A$ where $(N_\ga\colon \ga < \gl)$ is an increasing sequence of
 elementary submodels of cardinality less than $\gl$ of some  $({\rm
H}(\gk), {\in})$ 
with the usual good properties (including $X \in N_0$). Since 
the new filtration agrees with the old one on a closed unbounded set
we can assume $E$ consists only of ordinals of cofinality $\go$.

	Since the $p$-adic closure in $A$ of any set of cardinality
less than $\gl$ has cardinality less than $\gl$, $X$ has cardinality
$\gl$.  We can assume that the forcing is $\Fn(X, p, \go)$. Let $h$ be
the generic function from $X$ to $\oZ/p\oZ$. Note that the obvious
name for $h$ is in $N_0$. Suppose that $\hat h$ is a lifting of $h$.

	Choose a limit ordinal $\gd \in E$ and let $a$ be an element
of $A_{\gd + 1} \setminus A_\gd$ which is in the $p$-adic closure in
$A$ of $A_\gd$. Choose $q$ which determines $\hat{h}(a)$, say it
forces the value to be $m$. Choose $\ga < \gd$ so that $q \cap N_\gd =
q \cap N_\ga$. Call this restriction $q_1$.  Next take $n$ maximal so
that $a$ is $p^n$-divisible modulo $A_\ga$. (Since the $p$-adic
closure in $A$ of $A_\ga$ is contained in $A_\gd$ such an $n$ must
exist.) Choose $y \in A_\ga$ such that  $p^n \mid a - y$. Let $q_2 \in
N_\ga$ be an extension of $q_1$ which determines $\hat{h}(y)$. Now
choose $b \in A_\gd$ so that $p^{n+1} \mid a - y - b$. Notice that
$p^n$ is the exact power of $p$ which divides $b$ (even modulo
$A_\ga$).  Since $p^{n+1} \nmid b + A_\ga$, if we write $b/p^n$ as a
(infinite) sum of elements of $X$ there will be elements whose
coefficients are not divisible by $p$ which lie outside $A_\ga$. Hence
there is a condition $q_3 \in N_\gd$ extending $q_2$ so that $q_3$
forces $h(b/p^n) \not\equiv (m - \hat{h}(y)/p^n) \bmod p$. Let $r$ be
a common extension of $q_3$ and $q$. This gives a contradiction. 

To get the stronger statement on the cardinality of $\ext_p(A, \oZ)$ it is
enough as we did in Theorem~\ref{stdens-thm} to use the Cohen reals to code
$\gl$ generic functions and then note that for any two of them their
difference is also generic.
\fin

	Finally we can complete the proof of the main lemma.

\medskip
\proof (of Lemma~\ref{main}) The proof is by induction on $\gl < \gk$.
We show that any coseparable group $A$ of cardinality $\gl$ is free.
Suppose that $A$ is a coseparable group, $|A| = \gl < \gk$ and that
the induction hypothesis is true for all $\rho < \gl$. Assume, for the
sake of contradiction that $A$ is not free.  By the induction
hypothesis, $A$ is $\gl$-free and by the singular compactness theorem
\cite{sing} $\gl$ is regular. Let
$(A_\ga\colon \ga < \gl)$ be a $\gl$-filtration of $A$. There are two
main cases: either for a stationary set $E \se \gl$, $A_{\ga+1}/A_\ga$
contains a non-zero divisible group or not. In the first case we have
by Theorem~\ref{divis} that for every $p$ there is a counterexample to
the vanishing of $\ext_p(A, \oZ)$.

	The second case divides into two subcases. If there is a
stationary set $E$ so that for all $\ga \in E$, $A_{\ga+1}/A_\ga$ is
$\ha_1$-free then by the induction hypothesis and the fact that there
are only countably many primes we can choose a stationary subset $S$
of $E$ and a prime $p$ so that for all $\ga \in S$ there is a
counterexample to the vanishing of $\ext_p(A_{\ga+1}/A_\ga, \oZ)$. By
absorbing at most $\gl$ of the Cohen reals into the ground model we
can assume that $A$, $(A_\ga\colon \ga < \gl)$, $E$ and the
counterexamples to the vanishing of $\ext_p(A_{\ga+1}/A_\ga, \oZ)$ are
all in the ground model.  In this case we can apply
Theorem~\ref{nondivis}.

	In the final subcase
there is a stationary set $E$ so that for all $\ga \in E$,
$A_{\ga+1}/A_\ga$ is  not $\ha_1$-free and contains no non-zero
divisible subgroup. So for all $\ga \in E$ there is a prime $p_\ga$ so
that there is a counterexample to the vanishing of
$\ext_{p_\ga}(A_{\ga+1}/A_\ga, \oZ)$ (see \cite{EM} XII.2.7 for a
proof of this standard fact). If we choose a prime $p$ and $S$ a
stationary subset of $E$ so that $p_\ga =p$ for all $\ga \in S$, we
can again apply Theorem~\ref{nondivis} (as above). \fin

	Since we have been able to calculate the $p$-ranks of our groups in
all cases we have the following theorem on the structure of $\ext$. 

\begin{theorem}
Suppose it is consistent that a supercompact cardinal exists. Then it is
consistent with either $2^{\ha_0} = 2^{\ha_1}$ or $2^{\ha_0} < 2^{\ha_1}$
that for any group $A$ either $\ext(A, \oZ)$ is finite or $\ext(A, \oZ)$
has rank at least $2^{\ha_0}$.
\end{theorem}

	It is natural to ask whether the use of the large cardinals is
necessary. To show that it is, we will give a construction of a
coseparable group (which is similar to Chase's) from principles whose
negation implies consistency of the existence of many measurable
cardinals.  By the principle, $*(\gl, \ha_0)$, we will mean that there
exist a family $\set{S_i}{i<\gl^+}$ of countable subsets of $\gl$ such
that for any collection $I$ of size $\gl$ there exists $\set{S_i^*}{i
\in I}$ where for all $i \in I$ $S_i \setminus S_i^*$ is finite and
for all $i \neq j$, $S_i^* \cap S_j^* = \emptyset$.  (See page 157 of
\cite{EM} for details.)

\begin{theorem}
Suppose that $\gl$ is a cardinal, $\cf(\gl) = \go$, $2^\gl = \gl^+$,
$*(\gl, \ha_0)$ holds and  for all
$\mu < \gl$, $2^\mu < \gl$. Then there is a non-free
coseparable group of cardinality $\gl$. 
\end{theorem}

\proof Let $\set{S_i}{i < \gl^+}$ be as guaranteed by $*(\gl, \ha_0)$ and for
each $i$ let
$\set{s_n^i}{n < \go}$ be an enumeration of $S_i$.  Let $G_0$ be the group
freely generated by $\gl$ (= $\bigcup_{i<\gl^+} S_i$). Let
$\set{f_i}{i < \gl^+}$  enumerate  the homomorphisms from $G_0$ to
$\oZ/p\oZ$ as   
$p$ varies over the primes. We will inductively define a strictly
increasing chain 
of free groups $G_\ga$ and homomorphisms $g_\ga$ from $G_{\ga+1} \to \oZ$ for
$\ga < \gl^+$ so that $g_\ga$ is a lifting of $f_\ga$, $G_\ga/G_0$ is
divisible and for all $\gb < \ga$, $g_\gb$ extends to a homomorphism from
$G_\ga$ to $\oZ$. (Since $G_\ga/G_\gb$ is divisible there is at most one
extension.) If we can do this then  $G = \bigcup_{\ga <
\gl^+} G_\ga$ is the group required. $G$ is not free. Since $G/G_0$ is
divisible, 
for every prime $p$ every homomorphism from $G$ to $\oZ/p\oZ$ is uniquely
determined by its restriction to $G_0$. By the construction every 
homomorphism from $G_0$ to $\oZ/p\oZ$ lifts, so $G$ will be as required. 

	It remains to do the construction. We will work inside the
$\oZ$-adic completion  of $G_0$ so our use of infinite sums will be justified.
Choose a function $h\colon \gl \to \gl$ so that the inverse image of any
$\nu$ has cardinality $\gl$. 
For each $\ga$ we will choose $x^\ga_{n,0} \neq x^\ga_{n,1}$ so that
$h(x^\ga_{n,l}) = s_n^\ga$  and define
$x_\ga = 
\sum_{n<\go} n!(x_{n,0}^{\ga} - x_{n,1}^{\ga})$. The group $G_{\ga+1}$ will
be the pure closure of $G_\ga \cup \{x_\ga\}$. At limit ordinals we will
take unions. If we do this then for each $\ga$ the group $G_\ga$ will be
free but the group $G$ will not be free (see \cite{MS1} or \cite{EM} 
Theorem VII.2.13).

	Now we make the choices. Suppose we are at stage $\ga$. We have
already chosen $\set{g_\gb}{\gb < \ga}$ a set of at most $\gl$ functions.
Write $\ga = \bigcup_{n<\go} w_n$, where the cardinality of each $w_n$ is
less than $\gl$ and $w_n \se w_{n+1}$. Since $2^{|w_n|} < \gl$, for each $n$ we can choose
$x^\ga_{n,0} \neq x^\ga_{n,1}$ so that $h(x^\ga_{n,l}) = s_n^\ga$ and for all
$\gb \in w_n$, $ g_\gb(x^\ga_{n,0}) = g_\gb(x^\ga_{n,1})$. Then the
homomorphisms extend to the infinite sum (namely $x_\ga$) since each
homomorphism is 
$0$ on all but a finite number of the terms being summed. Finally we must
choose $g_\ga$. Since $G_{\ga+1}$ is free there is some lifting of $f_\ga$
to $G_{\ga+1}$, so we can choose $g_\ga$ as required. \fin

	As we have mentioned the statement that the hypothesis of the
theorem fails at all cardinals implies the consistency of many large
cardinals. So the use of some large cardinal hypothesis was necessary. One
might also want to know if the continuum needs to be as  large as we have
made it. The answer to this question is no. Assuming the consistency of the
existence of a supercompact cardinal it is possible to show that it is
consistent with $2^\gl = \gl^+ + \ha_2$ that every coseparable group is
free. To prove this result we will divide our effort between a
set-theoretic statement and a group theoretic one. First we need to
define a forcing axiom. By \ax\ we mean the statement ``if $\oQ$
is an $\ha_1$-complete forcing and $\tilde{S}$ is a $\oQ$-name and
$\force_\oQ\mbox{``} \tilde{S} \mbox{ is a stationary subset of
$\go_1$''}$ then for some 
directed $G \se \oQ$, $\set{\ga < \go_1}{\mbox{there is }q \in G,\
q\force \ga \in \tilde{S}}$ is stationary. ''

	The point of \ax\ is to give a  reflection
principle which shows that certain groups are not $\ha_2$-free.

\begin{lemma}
\label{reflection} Suppose that \ax\ holds and $\oP$ is an
$\ha_2$-complete forcing.  (By $\ha_2$-complete we mean that any
directed subset of cardinality $\ha_1$ has an upper bound.) If $G$ is
$\oP$-generic then the following statement is true in the generic
extension.
\begin{quote}
     Suppose $\gl$ is a regular cardinal, $A$ is a group of
cardinality $\gl$ and $\set{A_\ga}{\ga <\gl}$ is a $\gl$-filtration of
$A$ such that  
$$E = \set{\ga}{A_{\ga+1}/A_\ga \mbox{ is not $\ha_1$-free}}$$ 
is stationary. Then $A$ is not $\ha_2$-free.
\end{quote}
\end{lemma}

\proof Without loss of generality we can assume that the set
underlying $A$ is $\gl$ and $\tilde{A}$ is a $\oP$-name for  the
group structure on $\gl$ which is forced by the empty condition to be as
in the hypothesis. We will show that there is a condition which
forces that $\tA$ is not $\ha_2$-free. Let $\oQ$ be the product of
$\oP$ with the forcing $\oR$ which adds a generic function $\tf$ from
$\go_1$ onto $\gl$ by countable conditions. Let $\tilde{S}$ be the
name of a subset of $\go_1$ such that $q \force \ga \in \tS$ if and
only if for some $\gb > \ga$, $q$ determines $\tf\rest \gb$ and the
following holds. Let $X$
be  the determined value of the image of $\tf\rest \gb$. Then $q$
determines the group structure of $\tA$ on $X$ and if $Y$ is the
determined image of $\tf\rest \ga$ then $X$ and $Y$ are given group
structures by $q$ and $X/Y$ contains a non-free group of finite rank.

	We need to see that $\tS$ is forced to be stationary. First
recall that for any groups $C$ and $D$, where $C \se D$, $D/C$ is not
$\ha_1$-free if and only if there is a non-free subgroup of $D/C$ of
finite rank. Let $\tT$ be the $\oR$-name for $\tS^{G\times
\tilde{H}}$, where $\tilde{H}$ is the canonical name for an
$\oR$-generic set. It is enough
to see that in $V[G]$, $\tT$ is forced to be stationary. This task is
quite easy. Let $\tC$ be an $\oR$-name in $V[G]$ for a closed
unbounded subset of $\go_1$. By taking a new filtration of $A$ we can
assume that $E$ consists entirely of ordinals of cofinality $\go$. Let
$r$ be an element of $\oR$, i.e., $r$ is a function from some
countable ordinal $\mu$ to $\gl$. Choose $N$ an elementary submodel of
some appropriate $({\rm H}(\gk), \in)$ so that $|N| < \gl$ and $N \cap
A = A_\ga$ for some $\ga \in E$ also we can assume that $r$ and
everything else we have been talking about are elements of $N$. Choose
a countable subgroup $X_0$ of $A_{\ga+1}$ so that $X_0 + A_\ga /A_\ga$
is not free. It is standard to choose an increasing sequence $r = r_0,
r_1, \ldots$ of conditions in $N \cap \oR$ where the domain of $r_n$
is $\mu_n$ so that there is an increasing sequence $(\nu_n\colon n <
\go)$ of ordinals such that for all $n$, $r_n \force \nu_n \in \tC$,
$\mu_n < \nu_{n+1}$ and $\bigcup_{n < \go} {\rm rge} (r_n)$ is a
subgroup of $A_\ga$ which contains $X_0 \cap A_\ga$. Let $r_\go$
denote $\bigcup_{n < \go} r_n$.  Finally let $r_{\go+1}$ be an
extension of $r_\go$ so that the range of $r_{\go +1}$ is $X_0 + {\rm
rge}(r_\go)$.  Then $r_{\go+1}$ forces that $\bigcup_{n<\go} \nu_n \in
\tC \cap \tT$.

 So by \ax\
there is a directed set $H$ of cardinality $\ha_1$ such that
$\set{\ga}{\mbox{there is }q \in H, q \force \ga \in \tS}$ is
stationary. Let $p$ be an upper bound to the first coordinates of
$\{q\colon q \in H\}$. By the definition of $\tS$,  $H$ determines a
function, $f$ from $\go_1$ to $\gl$ and  a
group structure on the range of $f$ which $p$ forces to coincide with
a subgroup of $\tA$. Further the definition of $\tS$ and the choice of
$H$ guarantees that the group structure on the range of $f$ is not
free. \fin

\begin{lemma}
\label{comppres}
    Suppose that $\gl$ is a regular cardinal, $A$ is not coseparable
and $|A| = \gl$. If $\oQ$ is a $\gl$-complete notion of forcing then
$\oQ \force A \mbox{ is not coseparable}$.
\end{lemma}

\proof Without loss of generality we can assume that the set
underlying $A$ is $\gl$. Suppose $p$ is a prime and $h\colon A \to
\oZ/p\oZ$ is a 
function in the ground model which does not lift. Suppose that
$\tilde{g}$ is forced to be a lifting of $h$. Choose an increasing
sequence of conditions $(q_\ga\colon \ga < \gl)$ so that $q_\ga$
determines $\tilde{g}(\ga)$. Then $f$ defined by $f(\ga) = n$ if
$q_\ga \force \tilde{g}(\ga) = n$ is a lifting of $h$. \fin

\begin{theorem} Suppose that the following statements are true: for
all infinite cardinals $\gl$, $2^\gl = \gl^+ + \ha_2$; every
coseparable group of power less than $2^{\ha_0}$ is free; and \ax.
Then there is a generic extension such that 
\begin{trivlist}
\item[$(1)$]  cardinalities and cofinalities are preserved,
\item[$(2)$] for all infinite cardinals $\gl$, $2^\gl =
\gl^+ + \ha_2$, 
\item[$(3)$] every coseparable group is free, 
\item[$(4)$] suppose $\gl$ is a regular cardinal, $A$ is a group of
cardinality $\gl$ and $\set{A_\ga}{\ga <\gl}$ is $\gl$-filtration of
$A$ such that  $E = \set{\ga}{A_{\ga+1}/A_\ga \mbox{ is not
$\ha_1$-free}}$ is stationary. Then $A$ is not $\ha_2$-free. 
\end{trivlist} 
\end{theorem}

\proof The forcing is an iteration over all regular cardinals greater
than or equal $\ha_2$, where the support is an initial segment. The
forcing $\tilde{Q}_\gl$ is the $\oP_\gl$-name for adding a Cohen
subset of $\gl$. The result of this forcing we will use is that
$\dmd(S)$ holds for all regular $\gl > \ha_1$ and $S$ a stationary
subset of $\gl$. Properties (1) and (2) are standard. Property (4)
has been proved in Lemma~\ref{reflection}. It remains to verify that
every coseparable group is free. Since the forcing is
$\go_2$-complete, we have, by Lemma~\ref{comppres}, that every
coseparable group of cardinality at most $\ha_1$ is free. Suppose
that $A$ is a non-free coseparable group of cardinality $\gl$ where
$\gl$ is minimal. (Note that $\gl> \ha_1$.) We can assume that the
set underlying $A$ is $\gl$. 
By forcing with the initial segment which adds the subsets to $\mu$
for $\mu \leq \gl$, we can assume
we are working in a universe satisfying (1), (2), (4) where every
coseparable group of cardinality less than $\gl$ is free, $A$ is a
non-free $\gl$-free group and $\dmd(S)$ holds for every stationary
subset of $\gl$. 
Choose a $\gl$-filtration $(A_\ga\colon \ga < \gl)$ of $A$. The
filtration should be chosen so that $A_{\ga+1}/A_\ga$ is not free if
$A/A_\ga$ is not $\gl$-free. 

	Let $E = \set{\ga}{A_{\ga+1}/A_\ga \mbox{ is not
$\ha_1$-free}}$. Then $E$ is not stationary, since by (4), if $E$
were stationary $A$ would not be $\ha_2$-free. By the inductive
hypothesis, $\set{\ga}{\ext_p(A_{\ga+1}/A_\ga, \oZ) \neq 0}$ is
stationary for some $p$. (In fact, for all $p$.) Hence by
Theorem~\ref{nondivis} and the fact that $\dmd(S)$ holds for all
stationary subsets of $\gl$ we are done. \fin

	It remains to show the hypothesis is consistent. 

\begin{theorem}
Suppose it is consistent that a supercompact cardinal exists then it
is consistent that:  for
all infinite cardinals $\gl$, $2^\gl = \gl^+ + \ha_2$; every
coseparable group of power less than $2^{\ha_0}$ is free; and \ax.
\end{theorem}

\proof Since for any $\ha_1$-complete forcing $\oP$, there is some
cardinal $\gl$ so that $\oP \times \Fn(\go_1, \gl, \go_1)$ is
equivalent to $\Fn(\go_1, \gl, \go_1)$ it is enough to consider posets
of the form $\Fn(\go_1, \gl, \go_1)$ in the verification of \ax.
We will call $\Fn(\go_1, \gl, \go_1)$ the forcing for collapsing $\gl$
to $\go_1$.
Suppose that $\gk$ is supercompact. We can assume that GCH holds in
the ground model. Let $f:\gk \to \gk$ be a function
such that for any $\mu$ and $\gl$ there is some $M$ and an elementary
embedding of $j\colon V \to M$ so that $j(f)(\gk) = \mu$ and ${}^\gl M
\se M$. The forcing is an iteration $(\oP_i, \tQ_i\colon i < \gk)$,
where if $i = 2j$, $\tQ_i$ is the $\oP_i$-name for the forcing
collapsing $f(j)$ to $\go_1$ and for $i = 2j +1$, $\oQ_i$ is the
forcing for adding $\ha_1$ Cohen reals. A function $p$ is in $\oP_\gk$
if for all $i$, $p\rest i \force p(i) \in \tQ_i$, the support of $p$
intersect the even ordinals is countable and the support of
$p$ intersect the odd ordinals is finite.

It is standard to verify that
\ax\ is forced to hold. As well $\gk$ is forced to be either $\ha_1$
or $\ha_2$. We will show that $\ha_1$ is preserved by this forcing and
that every coseparable group of cardinality $\go_1$ is free in the
forcing extension. (Notice that the first statement is implied by the
second since otherwise CH would hold and so there would be a non-free
coseparable group of cardinality $\ha_1$. In fact we will verify that
$\go_1$ is preserved as part of the verification that all coseparable
groups of cardinality $\ha_1$ are free.) We will only sketch the proof
and leave details to the reader. We can reconstruct the iteration
differently. For any $i$, let $\oR_i$ be the finite support iteration
of $(\tQ_{2j+1}\colon 2j+1 < i)$. The forcing $\oP_\gk$ can be
rewritten as $\oR_\gk * \tS_\gk$, where $\tS$ is the countable support
iteration of $(\tQ_{2j}\colon 2j < i)$. We claim that if $H$ is an
$\oR_\gk$-generic set and $G$ is $\oP_\gk$-generic then $V[H]$ and
$V[G]$ have the same countable sets of ordinals. To show this we show
by induction on $i$, that  $V[H_i]$ and
$V[G_i]$ have the same countable sets of ordinals. Here $H_i$
(respectively,  $G_i$)
denotes the restriction of $H$ (respectively, $G$) to $\oR_i$
(respectively, $\oP_i$). In particular we have from the  induction
hypothesis that for any $\gl$, $\Fn(\go_1, \gl, \go_1)^{V[H_i]} =
\Fn(\go_1, \gl, \go_1)^{V[G_i]}$. So we can in the forcing replace
$\tQ_{2j}$ by the $\oR_{2j}$-name for the forcing for collapsing
$f(j)$ to $\go_1$, call this name $\tQ_{2j}'$. 

	Any element of ${\tQ_{2j}'}{}^{H_{2j}}$ of the form
$\set{\{(\check{\ga}, \check{\gb})\}\times A_{\ga\gb}}{\ga < \grg, \gb
< f(j)}$ where $\grg < \go_1$, for all $\ga$ and $\gb$ $A_{\ga\gb} \se
\oR_{2j}$, for all $\ga, \gb, \gt$, if $p \in A_{\ga\gb}$ and $q \in
A_{\ga\gt}$ then $p$ and $q$ are incompatible, and for all $\ga$,
$\bigcup{\gb<f(j)} A_{\ga\gb}$ contains a maximal antichain. Call such
names {\em pleasant\/}.  Notice that any pleasant name is forced by
the empty condition to be a function with countable domain from
$\go_1$ to $\gl$. As well, if $q$ is any name which is forced by the
empty condition to be a model then there is a pleasant name which the
empty condition forces to extend $q$.  The important property that we
will use  is that the union of a countable chain of pleasant names
which is forced by the empty condition to be increasing in stength is
again a pleasant name. With this observation in hand, we can prove the
claim. There are two cases to consider, the successor case and the
limit case. We will just do the limit case as the successor one is
similar. Suppose $i$ is a limit ordinal and $\tilde{g}$ is a
$\oP_i$-name for a function from $\go$ to the ordinals. It is enough
to show that there is a name $\tilde{s}$ in $\tS_i$ so that the empty
condition in $\oR_i$ forces that $\tilde{s}$ determines the value of
$\tilde{g}$. 
 By the discussion there is a
sequence $(\tilde{s}_n\colon n < \go)$ such that for all $n$ and $2j$
in the support of $\tilde{s}_n$, $\tilde{s}_n(2j)$ is a pleasant name,
the empty condition forces that  $\tilde{s}_n$ determines the value of
$\tilde{g}(n)$ and for all $n$ and $2j$ in the support of
$\tilde{s}_n$, the empty condition (in $\oR_{2j}$) forces that
$\tilde{s}_{n+1}(2j)$ extends $\tilde{s}_n(2j)$. Then $\tilde{s}$ can
be taken to be the coordinatewise union of the $\tilde{s}_n$.

	Similarly, we can prove the following fact.
\begin{quote}
    Suppose $\tilde{g}$ is an $\oP_\gk$-name for a function from
$\go_1$ to $\go_1$. Then there is a sequence $(\tilde{s}_\ga\colon
\ga < \go_1)$ of elements of $\tS_\gk$, so that the empty condition
in $\oR_\gk$ forces that the sequence is increasing and that
$\tilde{s}_\ga$ determines $\tilde{g}$ up to $\ga$.
\end{quote}
With this fact in hand we can prove that all coseparable groups of
cardinality $\ha_1$ are free. Suppose $A$ is a non-free coseparable
group in the generic extension. By doing an initial segment of the
forcing we can assume that $A$ is in the ground model. By
Lemma~\ref{main}, we know that $A$ is not coseparable in the generic
extension by $\oR_\gk$. Let $h:A \to \oZ/p\oZ$ be a function which
does not lift. By adding some of the Cohen reals to the ground model
we can assume that $h$ is in the ground model. Without loss of
generality, we can assume that the set underlying $A$ is $\go_1$.
Suppose that $\tilde{g}$ is forced to be a lifting of $h$ to $\oZ$.
Let $(\tilde{s}_\ga\colon\ga < \go_1)$ be as above. Then if $H$ is
$\oR_\gk$-generic, $H$ together with the sequence determines a
function, $\gf$ from $A$ to $\oZ$. Since the value of any 3 elements is
determined by some condition, this function must be a homomorphism which
lifts $h$. But $\gf$ is in $V[H]$ which contradicts the choice of $h$.

\section{$2^{\ha_0}$-free may imply free}

	In \cite{BD} it is shown that if $\gk$ is a supercompact
cardinal and $\gk$ Cohen reals are added to the universe then
$2^{\ha_0}$-free implies free. In order to prove Theorem~\ref{weak}, we
need to know that if we begin with a supercompact  cardinal $\gk$ and
first add $\gk^+$ Cohen subsets of $\go_1$ and then add $\gk$ subsets of
$\go$, then $2^{\ha_0}$-free implies free. The proof is relatively
standard but we will present it below since we need the corollary. 

\begin{lemma}
Suppose $A$ is a non-free group. Then in any extension of the universe
which preserves cofinalities and stationary sets, $A$ is non-free.
\end{lemma}

\proof Consider such an extension and suppose that $A$ is any non-free
group. By restricting to a non-free subgroup of $A$ we can assume that
$A$ is $\gl$-free and the cardinality of $A$ is $\gl$ or $A$ is
countable.  By the singular compactness theorem $\gl$ is regular. The
theorem is by induction on $\gl$. By Pontryagin's criterion, being
non-free is absolute upwards for countable groups. (More generally since
satisfaction of sentences in $L_{\go_1\go}$ is absolute being non-free
is absolute upwards for countable structures in any variety in a
countable language.) 

	Suppose we know the result for all $\gk < \gl$ and $A$ is
$\gl$-free, non-free and of cardinality $\gl$. Then there is a
$\gl$-filtration $(A_\ga\colon \ga < \gl)$ so that for a stationary set $E$,
$A_{\ga+1}/A_\ga$ is not free for all $\ga \in E$. By induction we know
in the extension $A_{\ga+1}/A_\ga$ is not free. Hence $A$ is not free
in the extension. \fin

	The lemma above applies in the more general situation of a
variety in a countable language.

\begin{lemma} 
Suppose that GCH is true in the ground model and let $\oP =
\Fn(I,2,\go_1) \times \Fn(J,2, \go)$.  Suppose $\oQ$ is $\Fn(\mu,2, \go_1)
\times \Fn(\rho, 2, \go)$ where $\mu, \rho > \ha_1$. If $G$ is $\oQ$
generic and  $K = V[G]$, then
in $K$, $\oP$ is $\ha_2$-c.c.\  and preserves stationary sets. 
\end{lemma}

\proof It is easy to show that $\oP$ is $\ha_2$-c.c.\ in $K$, since
$\oQ \times \oP$ is  $\ha_2$-c.c.\ over the ground model.  For
amusement, we will give the proof only with the assumption that $K$ is
an extension by an $\ha_2$-c.c.\ notion of forcing. Since $K$ is an
extension by an $\ha_2$-c.c.\ poset, for any large 
enough regular $\gk$ there is a club $C$ of subsets of ${\rm H}(\gk)$
such that for all $N \in C$, $G$ is generic over $N$. Suppose now that
$A$ is a maximal antichain in $\oP$ (in $K$) and $\tA$ is a name for
$A$.  Choose $N \prec {\rm H}(\gk)$ of cardinality $\ha_1$ so that $N$
is closed under sequences of length $\go$, $\tA \in N$ and $G$ is
generic over $N$. We will show that any element of $\oP$ is compatible
with an element of $N[G] \cap A$. Since $N[G]$ has cardinality
$\ha_1$, this suffices.  Consider any element $(p_1, p_2) \in
\oP$. Since $N$ is closed under countable sequences, $(p_1 \rest (I \cap
N), p_2 \rest (J \cap N)) = q$ is in $N$. Hence $q$ is compatible with
some element $r$ of $N[G] \cap A$. Finally $r$ and $(p_1, p_2)$ are
compatible.

	All that remains to see is that $\oP$ doesn't destroy any
stationary subsets of $\go_1$ (in $K$). Let $S$ be any
subset of $\go_1$ in $K$ which is forced by some condition in $\oP$ to
be non-stationary.  Since $\oQ$ is $\ha_2$-c.c.\  we can
choose a name $\tilde{S}$ for $S$ which uses only $\Fn((L, 2, \go_1)
\times \Fn(T,2, \go))$ where the cardinality of $L, T$ is $\ha_1$. 
Finally the homogeneity of $\Fn((L, 2, \go_1) \times \Fn(T,2, \go))\times
\Fn((I, 2, \go_1)\times \Fn(J,2, \go))$ implies that $S$ is not stationary
in $K$. \fin
 	
\begin{theorem} 
Suppose $\gk$ is a supercompact cardinal, $V$ satisfies
CH and $\oP = \Fn(\mu, 2, \go_1) \times \Fn(\rho, 2, \go)$, where $\mu,
\rho > \ha_1$.  Then $\oP$ forces that every $\gk$-free group is free. 
\end{theorem}

\proof Suppose that $G$ is $\oP$-generic and $A$ in $V[G]$ is
$\gk$-free.  Let $\gl = |A|$.  Let $\tA$ be a name for $A$.  Choose an
embedding $j:V \to M$ so that $j(\gk) > \gl$ and $M$ is closed
under sequences of length $\max \{\gk, \mu, \gl, \rho\}$. Let $H$ be a
$j(\oP)$-generic set which contains $j\mbox{''}G$. Notice that $A$ is
isomorphic to the interpretation of $j\mbox{''}\tA$ in $M[j\mbox{''}G]$
and $M[H]$. (We denote this image as $j\mbox{''}A$.) If $A$ is not free 
then $j\mbox{''}A$ is not free in $M[j\mbox{''}G]$ and hence by the two
lemmas above also not in $M[H]$. However $j$ extends to an elementary
embedding of $V[G]$ into $M[H]$. So $j(A)$ is $j(\gk)$-free and hence in
$M[H]$, $j\mbox{''}A$ is free. \fin

\begin{corollary} 
\label{free}
If it is consistent with ZFC that a supercompact cardinal exists then
both of the statements ``every $2^{\ha_0}$-free group is free and
$2^{\ha_0} < 2^{\ha_1}$'' and ``every $2^{\ha_0}$-free group is free and
$2^{\ha_0} = 2^{\ha_1}$'' are consistent with ZFC.  Furthermore, if it
is consistent that there is a supercompact cardinal then it is
consistent that there is a cardinal $\gk < 2^{\ha_1}$ so that if $\gk$
Cohen reals are added to the universe then every $\gk$-free group is
free. 
\end{corollary}
 

\end{document}